\documentclass[amsthm]{elsart}

\usepackage{amsmath, amssymb, amsfonts, amscd}

\newcommand{\Z}{\mathbb{Z}}
\newcommand{\C}{\mathbb{C}}
\newcommand{\R}{\mathbb{R}}

\newcommand{\F}{\mathbb{F}}

\newcommand{\cat}{\mathfrak{C}}
\newcommand{\D}{\mathfrak{D}}

\newcommand{\A}{\mathcal A_p}
\newcommand{\U}{\mathfrak{U}}
\newcommand{\UE}{\mathfrak{U}^{ev}}
\newcommand{\UB}{\mathfrak{U}/\beta}
\newcommand{\nil}{\mathfrak{Nil}}

\newcommand{\ob}{\mathcal{O}}
\newcommand{\tob}{\tilde\ob}

\newcommand{\RE}{R_1^{ev}}
\newcommand{\st}{St_1^{ev}}

\newcommand{\MUP}[1]{MU^*(#1)\hat{\otimes}_{MU^*} \F_p}

\newcommand{\p}[1]{{#1}^\wedge_p}

\newcommand{\CG}{\mathcal{C}(G)}
\newcommand{\CdG}{\mathcal{C}'(G)}
\newcommand{\GC}{G_{\mathbb{C}}}
\newcommand{\tG}{\tilde G}
\newcommand{\tGC}{\tG_{\C}}

\newcommand{\hl}{\bf}

\newcommand{\invlim}{\varprojlim}

\swapnumbers

\theoremstyle{plain}
\newtheorem{lemm}{Lemma}[section]
\newtheorem*{lem2}{Lemma}

\newtheorem{theom}[lemm]{Theorem}
\newtheorem*{thm2}{Theorem}

\newtheorem{propo}[lemm]{Proposition}
\newtheorem*{prop2}{Proposition}

\newtheorem{coro}[lemm]{Corollary}
\newtheorem*{coro2}{Corollary}

\theoremstyle{remark}
\newtheorem*{note2}{Note}

\theoremstyle{definition}
\newtheorem{rmk}[lemm]{Remark}
\newtheorem*{rmk2}{Remark}
\newtheorem{blank}[lemm]{}
\newtheorem{ex}[lemm]{Example}
\newtheorem*{ex2}{Example}

\begin{document}

\begin{frontmatter}

\title{Steenrod operations on the Chow ring of a classifying space}

\author{Pierre Guillot}

\address{Department of Pure Mathematics and Mathematical Statistics
       \\Cambridge CB30WB
       \\{\tt p.guillot@dpmms.cam.ac.uk}}

\begin{abstract} 
We use the Steenrod algebra to study the Chow ring $CH^*BG$ of the classifying space of an algebraic group $G$. We describe a localization property which relates a given $G$ to its elementary abelian subgroups, and we study a number of particular cases, namely symmetric groups and Chevalley groups. It turns out that the Chow rings of these groups are completely determined by the abelian subgroups and their fusion.
\end{abstract}

\end{frontmatter}

\section*{Introduction}

\begin{blank}Studying the ring $CH^*BG$ (see \cite{totaro}, \cite{pedro}) for a complex algebraic group $G$, leads one to the ill-formulated conclusion that the Chow ring is the ``nicest possible'' subring of the cohomology -- it is very often polynomial when $H^*BG$ has a lot of extra nilpotent elements. The action of the Steenrod algebra on these objects can be used to make this idea more precise.\end{blank}

\begin{blank}Let $\U$ be the category of unstable modules over $\A$, the mod $p$ Steenrod algebra, and let $\UE$ be the category of unstable modules which are concentrated in even dimensions. A typical object of $\U$ is thus the mod $p$ cohomology $H^*X$ of a topological space, while a typical object of $\UE$ is the mod $p$ Chow ring $CH^*X$ of a smooth algebraic variety $X$. It is well-known that the study of $\UE$ (at any prime) is similar in many ways to that of $\U$ at the prime 2. Indeed, the proofs of certain results in $\U$ at $p=2$ can be formally transposed to $\UE$ with practically no work, giving theorems for free. For example this is used extensively in \cite{lionel}. In a sense the present paper is another illustration of this, with a connection to algebraic geometry (and Chow rings in particular) which is probably new.

Since Lannes introduced the $T$-functor, it became clear that the cohomology ring $H^*\Z/p$ of a cyclic group of order $p$ played a crucial role in $\U$. In $\UE$, it is replaced by its even part, which is none other than $CH^*B\Z/p$. This is one reason why it appeared to us reasonable to suppose that Chow rings had some importance in $\UE$. The idea was that there should be some phenomena happening specifically in $\U$ at the prime 2 which should hold true in $\UE$ regardless of the prime if one only replaced the cohomology of elementary abelian $p$-groups, which are almost invariably involved in the statement of interesting results about $\U$, by the corresponding Chow rings.

Such an example indeed presented himself, which was the starting point of this work on classifying spaces. Namely, if $S_n$ denotes the symmetric group on $n$ letters, it is established in \cite{guna} that the decomposition
$$H^*BS_n=\invlim H^*BE$$ where $E$ runs through the elementary abelian $p$-subgroups of $S_n$, is valid exactly at the prime 2 and not for any other $p$. One of our results in the present paper is that $$CH^*BS_n=\invlim CH^*BE$$ for any prime -- see \ref{mainsym}.\end{blank}

\begin{blank}Since the existence of a decomposition as above may be elegantly regarded as a localization result in $\U$ or $\UE$, as we shall explain soon, the ``niceness'' of Chow rings alluded to in the opening paragraph can be expressed by saying that $CH^*BG$ is very often a ``local'' module in $\UE$ with respect to the subcategory of nilpotent modules, whereas $H^*BG$ is hardly ever local in $\U$ (all this will be made precise in the first section).

When $G$ is a group such that $CH^*BG$ is such a ``local'' module, there is surprisingly simple description of the Chow ring in terms of the cohomology ring: it will be written $$CH^*BG=L(\tob H^*BG)$$ where $\tob: \U\to \UE$ is right adjoint to the forgetful functor $\ob: \UE \to \U$ and $L: \UE\to\UE$ is a localization functor (details in \ref{chowfromcoho}).

This is not true of all groups $G$: the quaternion group is a counter-example (but see \ref{counter} where we argue that this might be specific to the prime 2). However, we have attempted to convince the reader that good examples of groups having this property abound. In the last two sections, we prove that most Chevalley groups (definition in \ref{mainchev}), including $GL_n$, $Sp_n$, $SO_n^+$, $SL_n$ when $p$ does not divide $n$, and all exceptional groups if $p>7$, as well as all products and certain quotients of these, are indeed local in the above sense. See section \ref{main}.\end{blank}

\begin{blank}In passing, note that  for $G$ a Chevalley group with maximal torus $T$ and Weyl group $W$, and for $k$ a finite field, things come down to proving that \begin{equation}\label{eq1} CH^*BG(k)=(CH^*BT(k))^W \end{equation} which parallels the result for the cohomology of compact Lie groups. The cohomology rings of the same finite groups certainly do not satisfy the analogous property.

That (\ref{eq1}) is equivalent to our localization property will follow from the group-theoretic result below, which had been obtained by Steinberg (\cite{steinberg}) but for which we give here a completely different, rather quick, proof which uses the computations in complex cobordism that we did in \cite{pedro}:

\begin{prop2} Let $G$ be a Chevalley group over a finite field $k$ of characteristic $\ne p$ containing the $p$-th roots of unity, and let $T$ be a split maximal torus. If we assume that each elementary abelian $p$-subgroup of $\GC$ (the associated Lie group over $\C$) is contained in a maximal torus, then each elementary abelian $p$-subgroup of $G(k)$ is conjugated to a subgroup of $T(k)$. 
\end{prop2}

We establish this in \ref{pedroborel}.\end{blank}

\begin{blank}The computations that we do here complete our previous results in \cite{pedro}. One general fact in this direction is that for a semi-simple Chevalley group $G$, the map $$CH^*BG(k)\to \MUP{BG(k)}$$ is always surjective, where $MU$ denotes complex cobordism. (The reader should note that the precise statements in sections 3 and 4 have reasonable restrictions on $p$ and $k$.) Combining this with our knowledge of the cobordism of Chevalley groups, we obtain that  the image of the cycle map $CH^*BG(k)\to H^*BG(k)$ coincides with that of $H^*B\GC \to H^*BG(k)$ for a certain ``Brauer lift'' $BG(k)\to B\GC$, where again $\GC$ is the associated complex Lie group over $\C$. Also, this image is a polynomial ring. See the proposition in section \ref{image}.\end{blank}

\begin{blank}{\hl Notations.}\label{intronotations} Let us have a discussion on Chow vs cohomology rings of abelian groups; it will mostly serve as a pretext to fix the notations for the rest of the paper, but might also help the reader understand what follows.

We will always use $H^*X$ to denote the mod $p$ cohomology of the space $X$; when $X=BG$ we might occasionally write $H^*G$. Similarly, $CH^*X$ denote the mod $p$ Chow ring of the smooth variety $X$ (=reduced scheme of finite type over $\C$). In \cite{fulton}, this would be $A^*X\otimes \F_p$. See \cite{totaro} for the definition of $CH^*BG$ obtained by approximating $BG$ by smooth algebraic varieties. We may here and there use the notation $CH^*G$. We point that, as we will see a module such as $CH^*X$ as an element in $\UE$, we will assign to an element in $CH^k X$ the degree $2k$, not $k$.

Let $p$ be odd. The cohomology of a cyclic group of order $p$ is $H^*B\Z/p=\F_p[v]\otimes\Lambda [u]$, while $CH^*B\Z/p=\F_p[v]$.
For an elementary abelian group we have similarly 
$$H^*(\Z/p)^n=\F_p[v_1,\cdots, v_n]\otimes \Lambda (u_1,\cdots, u_n)$$
and
$$CH^*(\Z/p)^n=\F_p[v_1,\cdots, v_n]$$

Let $\ob: \UE\to\U$ be the forgetful functor, and let $\tob$ denote its right adjoint; if $M$ is an unstable module, then $\ob\tob M$ is the largest submodule of $M$ which is concentrated in even dimensions. For example 
\begin{equation}\label{eq2} CH^*(\Z/p)^n=\tob H^*(\Z/p)^n \end{equation}

Now we make the following observation. The cohomology of a cyclic group $\Z/(p^r)$, for $r\ge 2$, has the same description as in the case $r=1$ above. However, as module over $\A$, there are differences: for example, for $\Z/p$ we have $\beta (u)=v$, where here and elsewhere $\beta$ denotes the Bockstein, whereas for $\Z/(p^2)$ we have $\beta (v)=0$. On the other hand, the Chow rings of these groups agree as objects of $\UE$. The relation (\ref{eq2}) above therefore does not hold for general abelian groups. (We will prove in this paper that it does hold for $S_n$. It is rather exceptional.)

Further, the restriction map induced by the inclusion $\Z/p\to \Z/(p^r)$ is not in general an isomorphism between the cohomology rings. However this holds true for Chow rings, and we can and will replace any abelian group by its maximal elementary abelian $p$-subgroup when it comes to computing the mod $p$ Chow ring.

At $p=2$ now, $H^*\Z/2$ is a polynomial ring on one variable $u$ in degree 1. We put $v=u^2=\beta(u)=Sq^1(u)$ to have uniform notations at all primes.

\end{blank}

\begin{blank}{\hl Organization of the paper}. In the first section we discuss localization in the categories $\U$ and $\UE$ in some detail, in particular we explain the relation between the localization of $H^*BG$ or $CH^*BG$ and the elementary abelian subgroups of $G$. In section 2 we study the particular case $G=S_n$. Section 3 prepares the ground for the study of Chevalley groups, discussing Weyl groups, universal covers, and reducing the localization of these groups to the statement already mentioned. It is in section 4 that explicit computations are made with them. There is also an appendix on the double coset formula, a handy technical tool used in several places throughout the article.
\end{blank}

\section{Nilpotent modules and localization}\label{section:nil}

\begin{blank}
Let $G$ be a reductive algebraic group, and let $\CG$ be the category whose objects are the elementary abelian $p$-subgroups of $G$ and whose morphisms are induced by conjugations in $G$. In \cite{quillen}, Quillen proved that the natural map
$$H^*(BG)\to \underset{\CG}{\invlim} H^*BE$$ is an ``F-isomorphism'', in the sense that any element in the kernel is nilpotent, and for any $x$ in the target, $x^{p^n}$ is in the image for some $n$. As it happens, this is a property that can be characterised in $\U$, because the $p$-th power is given by a Steenrod operation. Quillen's theorem can be reformulated by saying that the above map  becomes an isomorphism upon localizing away from the subcategory of so-called nilpotent modules (we will give more precise definitions below). Rather surprisingly, it is very fruitful to ``linearize'' the situation in this way (\cite{guna}, \cite{localize}). 

Yagita proved a version of Quillen's theorem for Chow rings when $G$ is finite, see \cite{yagita}: the map $$CH^*(BG)\to \underset{\CG}{\invlim} CH^*BE$$ is also an F-isomorphism. We will see in this section that we can express this as a localization result in $\UE$, parallel to the one for cohomology. If we saw our Chow rings as modules in $\U$ this would not be possible, which confirms our intuition that Chow rings naturally live in $\UE$.

It is natural to ask whether there are any groups for which Quillen's map is actually an isomorphism. For example, as already announced in the introduction, if $G=S_n$ is the symmetric group, the map is an isomorphism at the prime 2 but not for any other prime, see \cite{guna}. However we shall see in the next section that the Quillen-Yagita map for $S_n$ is an isomorphism at all $p$. Indeed the main point of this paper is to show that this is true for a whole lot of groups.

In this section we make the above statements precise, treating simultaneously the cases of $\U$ and $\UE$. There is essentially nothing new here, cf \cite{localize}, \cite{guna}, \cite{preprint}, but we need all this for reference. It is also felt that the reader might appreciate a concise and reasonably self-contained presentation.

\end{blank}

\begin{blank}{\hl Quotient categories; localizations.} (\cite{grot}, \cite{gabriel}) Let $\cat$ be an abelian category, and let $\D$ be a Serre class (a full subcategory with the property that if two objects of a short exact sequence in $\cat$ belong to $\D$, then so does the third). Then there is a quotient category $\cat/\D$ which is abelian and an exact functor $r: \cat\to\cat/\D$ satisfying the obvious universal property. A morphism in $\cat$ induces a monomorphism (resp. epimorphism) in $\cat/\D$ if and only if its kernel (resp. cokernel) belongs to $\D$. We will talk about $\D$-monomorphisms, etc.

The category $\D$ is said to be {\em localizing} if $r$ admits a right adjoint, or section functor, $s: \cat/\D\to \cat$. We put $l=s\circ r$ and we call the natural map $\lambda_M: M\to l(M)$ the {\em localization of $M$ away from $\D$}. When this is an isomorphism, we say that $M$ is $\D$-closed or $\D$-local.
 
One can prove easily that the natural transformation $r\circ s (M)\to M$ is an isomorphism. It follows that $\lambda_M: M\to l(M)$ is a $\D$-isomorphism, that is, $r(\lambda_M)$ is an isomorphism. Hence $\lambda_{l(M)}=l(\lambda_M): l(M)\to l\circ l (M)$ is an isomorphism, too: in other words, the localization of $M$ is local. More generally, keep in mind that a $\D$-isomorphism $M\to N$ induces an isomorphism $l(M)\to l(N)$.

Finally, we note that if $\cat$ has enough injectives, then an object $M$ is $\D$-closed if and only if $$Ext^i_\cat(D,M)=0$$ for $i=0,1$ and all $D$ in $\D$.

\end{blank}

\begin{blank}{\hl Nilpotent and $\nil$-closed modules.}\label{nildefs} Given $M$ in $\U$ or $\UE$, and $x\in M$ of even dimension, we put $P_0x=P^{|x|/2}x$ (so that if $M$ happens to be an unstable algebra, $P_0x=x^p$).

We say that $M$ is {\em nilpotent} if $P^N_0x=0$ for all $x$ of even dimension and all large $N$. (See \cite{lionel}, p47). There is an exception if we work with $\U$ at the prime 2: in this case we put $Sq_0x=Sq^{|x|}x$ and call a module $M$ nilpotent if $Sq^N_0x=0$ for all $x$ and all large $N$.

The subcategory of $\U$ or $\UE$ which is comprised of the nilpotent modules will be denoted by $\nil$ (we do not distinguish between $\U$ and $\UE$, hopefully the context will make things clear).

It is easy to see that $\nil$ is localizing, see the criterion in \cite{lionel}, prop 6.3.1. Therefore, the results of the previous paragraph apply, and we use the notation $\lambda_M: M\to L(M)$ for the localization away from $\nil$ (notation as in \cite{localize}).

It is well-known that $\U$ has enough injectives, and it follows that the same can be said of $\UE$. Accordingly, $M$ is $\nil$-closed if and only if $Ext^i(N,M)=0$ for $i=0,1$ and for all $N$ nilpotent.

To finish with, we call a module {\em reduced} if $Hom(N,M)=0$ for all nilpotent modules $N$. In $\UE$, or in $\U$ at the prime 2, this is equivalent to demanding that $P_0$ be injective on $M$: combine lemma 2.6.4 and equation 1.7.1* in \cite{lionel} (alternatively, see lemma 4.5 in \cite{preprint}). 

\end{blank}

\begin{lemm} Any reduced module in $\UE$ embeds in a reduced $\UE$-injective. The tensor product of two reduced $\UE$-injectives is a reduced $\UE$-injective.
\end{lemm}

\begin{proof} We observe immediately that if $I$ is injective in $\U$, then $\tob I$ is injective in $\UE$. The first assertion of the lemma follows then from the corresponding statement in $\U$, which is well-known, and the left exactness of $\tob$. It also follows that any reduced injective in $\UE$ is a direct summand in such a module $\tob I$ with $I$ reduced, and we get the second assertion, again because the analogous result in $\U$ is well-known (and because $\tob (A\otimes B)=\tob A \otimes \tob B$ if one of the factors is reduced).
\end{proof}

\begin{propo}\label{nilprop} In either $\U$ or $\UE$, we have the following properties:
\begin{enumerate}
\item If $0\to M' \to M \to M'' \to 0$ is exact and if $M'$ and $M''$ are both $\nil$-closed, so is $M$.
\item If $0 \to M \to M' \to M''$ is exact and if $M'$ is $\nil$-closed while $M''$ is reduced, then $M$ is $\nil$-closed.
\item $M$ is $\nil$-closed $\iff$ there exists an exact sequence $$0\to M \to I \to J$$ with $I$ and $J$ both reduced and injective.
\item $M_1$ and $M_2$ $\nil$-closed $\Rightarrow$ $M_1\otimes M_2$ $\nil$-closed.
\end{enumerate}

\end{propo}

\begin{coro}\label{invlim} Any product or inverse limit of $\nil$-closed modules is $\nil$-closed.
\end{coro}

\begin{proof} This is trivial from the lemma. Note that the third point in the proposition is only here to prove the fourth. In case of problems see \cite{guna} or \cite{preprint}.
\end{proof}

\begin{blank}{\hl Back to Quillen's map.} The cohomology of an elementary abelian $p$-group is reduced and injective (\cite{lionel}, 2.6.5 and 3.1.1), hence is $\nil$-closed. Consequently, the target of Quillen's map is $\nil$-closed, as an inverse limit of such. Quillen's theorem asserts that his map is a $\nil$-isomorphism, so that it becomes an isomorphism upon localizing. In other words $$L(H^*(BG))\approx\underset{\CG}{\invlim} H^*BE$$ (as the inverse limit is isomorphic to its own localization). We see in this way that Quillen's map is an isomorphism if and only if $H^*BG$ is $\nil$-closed.

Similarly, $CH^*BE=\tob H^*BE$ is injective and reduced in $\UE$, and the target of the Quillen-Yagita map is $\nil$-closed in $\UE$. It follows that  $$L(CH^*(BG))\approx\underset{\CG}{\invlim} CH^*BE$$ and that the Quillen-Yagita map is an isomorphism if and only if $CH^*BG$ is $\nil$-closed.

\end{blank}

\begin{ex} It is instructive to have a look at $G=GL(n,\C)$, even if this $G$ is not finite (so that the Quillen-Yagita map cannot be mentioned): in effect $M=H^*BG=CH^*BG$ illustrates well how a module can behave better in $\UE$ than in $\U$, which we see as follows. If $E$ is the subgroup of $G$ of diagonal matrices with $p$-th roots of unity as entries, then any elementary abelian subgroup of $G$ is conjugated to a subgroup of $E$. If $W$ denotes the Weyl group, it follows that $$L(M)=(H^*BE)^W\ne H^*BGL(n,\C);$$ consequently $M$ is not $\nil$-closed in $\U$. On the other hand $$(CH^*BE)^W=CH^*BGL(n,\C)$$ from which we deduce that $M$ is indeed $\nil$-closed in $\UE$ (being an inverse limit of $\nil$-closed modules).

To get an example with finite groups, it is in fact possible to take a finite field $k$ of characteristic different from $p$ but containing the $p$-th roots of unity, and to consider $H^*BGL(n,k)$ and $CH^*BGL(n,k)$. The former is reduced but not $\nil$-closed in $\U$ (\cite{localize}, \cite{quilgln}), while the latter is $\nil$-closed in $\UE$ (see \cite{pedro}). In the last section of this article we shall extend this to other groups of matrices over finite fields. 

\end{ex}

\begin{ex}\label{counter} We cannot hope for $CH^*BG$ to be always $\nil$-closed, as is illustrated by the group of quaternions: there is only one elementary $2$-subgroup, namely the centre $\{1;-1\}$, and the restriction map is neither injective nor surjective.

However, this example is not too discouraging: the quaternion groups are precisely the $2$-groups which have $2$-rank 1, together with cyclic groups (recall that the rank of a $p$-group is the dimension of the largest $\F_p$ vector space contained in it). For odd $p$ on the other hand, this pathology disappears, and a group of $p$-rank 1 can only be cyclic. So the example above might reflect a purely group-theoretic defect (or subtlety, if you want), and one might still hope that for groups of odd order, the Quillen-Yagita map is ``very often'' an isomorphism. 
\end{ex}

\begin{blank}{\hl Remark.}\label{chowfromcoho} Combining Quillen's and Yagita's result, we obtain of course that the cycle map $CH^*BG\to H^*BG$ is an F-isomorphism. So we can write in $\U$ that $L(CH^*BG)=L(H^*BG)$ and in $\UE$ that $L(CH^*BG)=L(\tob H^*BG)$. Whenever $CH^*BG$ is $\nil$-closed this reads $$CH^*BG=L(\tob H^*BG)$$
This is a description of the Chow ring in terms of $H^*BG$ using merely functors between $\U$ and $\UE$ (which do not depend on $G$).

We end this section with a few more remarks, before starting to give examples of groups $G$ with $CH^*BG$ $\nil$-closed, hopefully convincing the reader that there is a fair number of them. When $p$ is odd, we do not know of an example of a group {\em not} satisfying this property.

\end{blank}

\begin{blank}{\hl Remark.}\label{tob} Suppose that $G$ is a group such that $H^*BG$ is reduced. Since the functor $\tob$ is left-exact and commutes with inverse limits, we have an injection:
$$\tob H^*BG \hookrightarrow \invlim \tob H^*BE$$ But $\tob H^*BE=CH^*BE$, so if we suppose further that the Quillen-Yagita map for Chow rings is surjective, it follows that the image of $CH^*BG$ under the cycle map has to contain all of $\tob H^*BG$. We will use this later with $G=S_n$.
\end{blank}

\begin{blank}\label{cobordism}{\hl Complex cobordism.} Recall (\cite{totaro}) that the cycle map factors as
$$CH^*BG \to \MUP{BG} \to H^*BG$$
where $MU$ denotes complex cobordism. Most of what we have said so far applies to the theory $\MUP{-}$.

We start by explaining why we have Steenrod operations on this theory. First, we have the isomorphism {\em of graded abelian groups}:
$$MU^*(MU)=MU^*(pt)[[\cdots c_n \cdots]]$$ More generally for any oriented cohomology theory $E$, the Thom-Dold theorem asserts $E^*MU=E^*BU$, and $E^*BU$ is computed via the Atiyah-Hirzebruch spectral sequence, and is a power series ring as above. {\em However this is not a ring isomorphism:} we are interested in the multiplication given by composition in $MU^*(MU)=[MU,\Sigma^* MU]$. Nevertheless, it is a fact (cf \cite{rudyak},VII,3.1) that any element $x\in MU^*(MU)$ can be written $$x=\sum_\omega a_\omega S_\omega$$ where $a_\omega\in MU^*(pt)$ and $S_\omega$ has an (easy) explicit description in terms of the $c_i$'s. (The ring structure on $MU^*(pt)$ is the natural one, and the notation $a_\omega S_\omega$ does refer to the ``correct'' multiplicative structure.) From this it follows that the kernel $I$ of the (surjective) map $MU^*(MU)\to H^*MU$, where as always $H$ denotes mod $p$ cohomology, is the ideal generated by $(p,x_i)\subset MU^*(pt)$, in standard notation. In other words, $I$ is generated in $MU^*(MU)$ by the kernel of the map $$MU^*(pt)=\Z[x_1,x_2,\cdots]\to \F_p=H^*(pt)$$

Now, this $I$ acts trivially on $\MUP{X}$ for any space $X$, of course. (That is, the elements of this ideal act as $0$.) Thus $H^*MU$ acts on our theory (again, with some ring structure which is not the one coming from the isomorphism $H^*MU=H^*BU$, but this will not matter). Since the Steenrod algebra acts on $H^*MU$, it also acts on $\MUP{X}$. Moreover, it is well-known that $H^*MU$ is a free $\mathcal A_p /(\beta)$-module, so that the Bocksteins act trivially in this new setting. (Note that $H^*MU$ is thus {\em not} unstable, this is an example that shows that the Thom-Dold isomorphism is not a map of $\mathcal A_p$-modules!)

Is $\MUP{X}$ unstable? This appears to be true, and we outline a proof -- we shall not need the result in what follows, as we shall only consider some spaces $X$ for which $\MUP{X}$ injects into the cohomology ring. To prove the claim, one would construct directly some operations on $\MUP{-}$ in the obvious way: starting with a manifold $M$ representing a class in $MU^*X$, then $M^p$ gives a class in $MU^*X^p$, and thus it yields one in $MU^*Z^pX$, the cyclic product; project this onto $$\MUP{X\times B\Z/p}=\left( \MUP{X} \right) \otimes \left( \MUP{B\Z/p} \right) $$ and define the coefficients of the polynomial in $v$ thus obtained to be your ``Steenrod'' operations. These act ``unstably'' for any space, by construction. After taking limits, we extend this to $X=MU$, and the images of $1\in \MUP{MU}$ define elements $\bar P^i$ there; these lift to $MU^*(MU)$ and yield operations on $MU^*(-)$ and on $\MUP{-}$ which by naturality agree with the ones just described. Thanks to the existence of such lifts, in fact, any element in $\MUP{MU}$ acts on $\MUP{-}$, and thus any relation between the $\bar P^i$'s and the $c_i$'s still holds after it's applied to any $x$ in $\MUP{X}$ for any $X$. By picking some spaces $X$ for which $\MUP{X}$ injects into the cohomology $H^*X$, and using the obvious compatibility of the $\bar P^i$'s with the actual Steenrod operations $P^i$ once we have projected to cohomology, we deduce $\bar P^i=P^i$.

In any case, if we take this for granted, we see that the modules $\MUP{X}$ live naturally in the category $\UB$ of unstable modules over $\mathcal A_p/(\beta)$, which is an intermediate between $\UE$ and $\U$. Considering the odd and even parts of a module $M$ in $\UB$ gives a decomposition of $M$ as a direct sum of two elements of $\UE$. It should be easy from this to extend the localization results obtained for $\U$ and $\UE$ to the category $\UE$. However, we do not have any application in sight, mostly because all modules of interest to us given as $\MUP{BG}$ for some $G$ are concentrated in even dimensions and can be seen as objects of $\UE$. Therefore we shall not fill in the details here.

 In any case, it is convenient to use the term ``F-isomorphism'', and we state for future reference Yagita's result just alluded to:

\begin{prop2}
Let $G$ be a finite group. The natural map 
$$\MUP{BG}\to \invlim \MUP{BE}$$ is an F-isomorphism.
\end{prop2}

\begin{coro2} If $G$ is a finite group, then the two maps
$$CH^*BG\to \MUP{BG} \to H^*BG$$ are F-isomorphisms.
\end{coro2}

\end{blank}

\section{The Symmetric Groups}

\begin{blank}\label{intro} In this section we prove that $CH^*BS_n$ is $\nil$-closed in $\UE$. As one could expect, the proof is by induction, by proving that if $CH^*BG$ is $\nil$-closed, then so is $CH^*B(S_p\wr G)$. It might seem simpler to use $\Z/p\wr G$ instead, but we want the reader to be able to compare our proof with that in \cite{preprint} which deals with a similar result using $\tob H^*BG$ instead of our $CH^*BG$ (and in turn, both treatments follow closely the original one in \cite{guna} for mod 2 cohomology). Incidentally, the result in \cite{preprint} can be recovered from ours, see \ref{feedback}.

For technical reasons stemming from \cite{totaro}, we will have to restrict attention to a certain class of groups: namely, $G$ will be assumed to have a subgroup $H$ of index prime to $p$ such that $BH$ can be cut into open subsets of affine space. Then $\Z/p\wr H$ is a subgroup of index prime to $p$ in $S_p\wr G$ which has the same property, by \cite{totaro}, lemma 8.1, and the induction may proceed.

We will also assume that the mod $p$ cycle map $CH^*BG\to H^*BG$ is injective. We will prove quickly that the cycle map for $\Z/p\wr G$ and $S_p\wr G$ is injective too.

We begin by explaining how Chow rings and cohomology rings are affected by taking wreath products, and we give a few immediate properties, in particular we will end up with an exact sequence which together with \ref{nilprop},1/, will eventually yield the result.
\end{blank}

\begin{blank}\label{bigraded}{\hl Cohomology and Chow rings of cyclic products.} We shall need the following nice result of Nakaoka\cite{nakaoka}:
$$H^*(\Z/p\wr G)=H^*(\Z/p, (H^*G)^{\otimes p})$$ 

Recall that the cohomology of a cyclic group with coefficients in any ring $A$ is periodic, and the period is given by taking cup-products with an element in $H^2(\Z/p,A)$. Here we denote by $v\in H^2(\Z/p, (H^*G)^{\otimes p})$ an element giving the period (abusing the notation given in the introduction).

For Chow rings of cyclic products, we use results of \cite{totaro}. There a certain functor $F_p$ from graded abelian groups to graded abelian groups is defined, which comes equipped with a natural map $F_pCH^*X\to CH^*Z^pX$, where $X$ is a variety and $Z^pX$ is it $p$-fold cyclic product. After changing the grading from dimension to codimension, reducing modulo $p$, and changing the notations to relate to \cite{bros}, the definition of $F_p A^*$ is as follows: take the $p$-fold tensor product $A^*\otimes\cdots\otimes A^*$, take symbols $Px$ in degree $p\cdot |x|$ for $x\in A^*$ of positive degree, take also elements $\alpha_i x$ of degree $p\cdot |x| + 2i$ for all $x\in A^*$ and positive $i$, and finally divide by the relations:
$$x_1\otimes \cdots \otimes x_p = x_2\otimes\cdots\otimes x_p\otimes x_1$$
$$x^{\otimes p}=0$$
$$P(x + y)= Px + Py + \sum s_1\otimes \cdots \otimes s_p$$
together with the relations that turn $\alpha_i$ into a homomorphism of groups. Here the sum in the third formula runs over a set of representatives for the $\Z/p$-orbits in the set $\{x,y\}^p - \{(x,\cdots,x),(y,\cdots,y)\}$. To be rigorous, the elements $\alpha_i x$ should be only defined for $i$ small enough and likewise the $Px$'s should only appear when the degree of $x$ is small enough; but in both cases the bound depends on the dimension of the variety $X$, and when we deal with classifying spaces of groups $BG$ we take a limit of varieties (with ``compatible'' Chow rings) having their dimensions going to infinity, and we do not need to worry about this complication.

From the construction of $F_p$ as given in \cite{totaro}, it is easy to describe the composition $F_p CH^*BG \to CH^*B(\Z/p\wr G)\to H^*B(\Z/p\wr G)$, for any $G$: the element $x_1\otimes\cdots\otimes x_p$ is sent to the ``norm'' $\sum x_{i_1}\otimes\cdots\otimes x_{i_p}$, with the sum running over all cyclic permutations, sitting in $H^0(\Z/p,(H^*G)^{\otimes p})$; the element $Px$ is sent to $x\otimes\cdots\otimes x$ in the same group (note how the relations above relate to the expansion of $(x+y)\otimes\cdots\otimes (x+y)$); and finally the element $\alpha_i x$ goes to $v^i\cdot x\otimes\cdots\otimes x$. This is all clear. In particular this map is injective if the cycle map for $G$ is injective.

For $G$ as in the introduction, lemma 8.1 in \cite{totaro} implies that $F_pCH^*BG\to CH^*B(\Z/p\wr G)$ is surjective. Therefore for the $G$ we consider, there is an isomorphism $F_pCH^*BG=CH^*B(\Z/p\wr G)$, the cycle map for $\Z/p\wr G$ is injective, and its image is explicitly described.

We shall denote by $\tau$ the transfer from $G^p$ to $\Z/p\wr G$. Its image is spanned by the ``norms'', clearly.

\end{blank}

\begin{blank}

We will also be interested in wreath product of the form $S_p\wr G$. Nakaoka has established that $H^*(S_p\wr G)=H^*(S_p, (H^*G)^{\otimes p})$ in this case too. Since $\Z/p$ is a $p$-Sylow of $S_p$, we deduce:

\begin{lem2}\label{swan} Let $W=N_{S_p}(\Z/p)/\Z/p=(\Z/p)^*$. Then
$$H^*(S_p\wr G)=(H^*(\Z/p\wr G))^W$$
and
$$CH^*(S_p\wr G)=(CH^*(\Z/p\wr G))^W$$
\end{lem2}
\end{blank}

\begin{proof} The equality for cohomology groups follows from Nakaoka's results just quoted and Swan's lemma. To get the result for Chow rings, we proceed as follows: there is an obvious inclusion, and the double coset formula (see the appendix) tells us that the image of the restriction map is the group of ``stable'' elements, ie those $x$ such that, putting $K=\Z/p\wr G$, we have $x|_{K\cap gKg^{-1}}=gxg^{-1}|_{K\cap gKg^{-1}}$ (proposition \ref{stable}). So we need to show that if $x$ is $W$-invariant, then $x$ is stable. But $K\cap gKg^{-1}$ is either $K$ or $G^p$ (this is because $G^p$ is normal in $S_p\wr G$, and the order of $\Z/p$ is prime). In either case, the cycle map is injective on this group, and so the result for cohomology implies that for Chow rings.
\end{proof}

\begin{propo} There are exact sequences:
$$0\to \tau(H^*(BG^p))\to H^*B(\Z/p\wr G)\to H^*BG\otimes H^*B\Z/p$$
and
$$0\to \tau(CH^*(BG^p))\to CH^*B(\Z/p\wr G)\to CH^*BG\otimes CH^*B\Z/p$$
\end{propo}

Note that the maps on the right come from the inclusion of $\Z/p\times G$ in $\Z/p\wr G=\Z/p\ltimes G^p$.

\begin{proof} The exact sequence for cohomology follows from results of Steenrod's (\cite{steenrod} chapter VII), as is explained, for example, in \cite{mui}, theorem II.3.7.

To get the result for Chow rings, we use the cycle map, which here is injective. The only thing to prove is that $CH^*B(\Z/p\wr G)\cap \tau(H^*(BG^p))=\tau(CH^*(BG^p))$, but this is clear from the explicit description in \ref{bigraded}.
\end{proof}

\begin{coro} There are exact sequences:
$$0\to (\tau(H^*BG^p))^W\to H^*B(S_p\wr G)\to R_1(H^*BG)\to 0$$
and
$$0\to (\tau(CH^*BG^p))^W\to CH^*B(S_p\wr G)\to \RE(CH^*BG)\to 0$$
where $W$ is as in lemma \ref{swan} and where $R_1(H^*BG)$, resp. $\RE(CH^*BG)$, is a submodule of $H^*BG\otimes H^*B\Z/p$, resp. $CH^*BG\otimes CH^*B\Z/p$.
\end{coro}

The functorial notations $R_1(H^*BG)$ and $\RE(CH^*BG)$ will be justified below.

\begin{blank} The point now is to prove that $(\tau(CH^*BG^p))^W$ and $R_1^{ev}(CH^*BG)$ are both $\nil$-closed, and to use proposition \ref{nilprop}, 1/. This is where we can only continue the proof for Chow rings, as the result does not hold for cohomology at odd primes. It is proved in \cite{preprint}, however, that if one replaces $H^*BK$ by $\tob H^*BK$ for all groups $K$ occuring in Quillen's map, then one still gets an isomorphism.
\end{blank}

\begin{blank}{\hl The functor $\RE$.} Let $M$ be a module in $\UE$. For any $x\in M$ of degree $2k$, define $$\st(x)=\sum_{i=0}^k (-1)^i v^{i(p-1)}\otimes P^{k-i}x \in CH^*B\Z/p \otimes M$$ where we recall that $P^i=Sq^{2i}$ when $p=2$. We define $\RE M$ to be the $CH^*(BS_p)$-submodule of $CH^*B\Z/p\otimes M$ generated by the elements $\st x$, for $x\in M$. Here we view $CH^*B\Z/p$ as a module over $CH^*BS_p$ via the restriction map.

There is a functor $R_1: \U\to\U$ which is defined in a similar, but more complicated, way (see for example \cite{preprint}, definition 4.2, or the original in \cite{zarati}). In fact one has $\RE M = \tob R_1(\ob M)$: this can be seen from the explicit description of $R_1$ given in {\em loc cit} , and using lemma \ref{swan} which asserts in particular that $CH^*S_p=\F_p[v^{p-1}]$. This proves that $\RE M$ is always in $\UE$, ie it is stable under the action of the Steenrod algebra. We will not use this seriously, however, and $\RE$ can be taken as a functor from $\UE$ to graded $\F_p$-vector spaces.

The definition of $\RE$ is very explicit and will allow computation. However:

\begin{lemm} The two definitions of $\RE M$ given coincide.
\end{lemm}

\begin{proof} We need to show that the image of the map $$CH^*(S_p\wr G)\to CH^*B\Z/p\otimes CH^*G$$ is $\RE M$. The fact that the elements $\st x$ are in this image follows from the very definition of the Steenrod operations on Chow rings: in \cite{bros}, an element $P(x)$ is constructed in $CH^*(S_p\wr G)$ (prop 4.2 in {\em loc cit}) which retricts to $\st x$ (definition 7.5 in {\em loc cit}; note that our $(-1)^i$ sign is on p10, before prop 6.6 there). Therefore, the image certainly contains $\RE M$.

To get the reverse inclusion, observe that (with notations as in \ref{bigraded}, our $Px$ being consistent with Brosnan's) the elements in $CH^*(\Z/p\wr G)$ not mapping to 0 in $CH^*G\otimes CH^*B\Z/p$ are of the form $Px$ or $\alpha_ix=v^i\cdot Px$. The latter elements can only be in $CH^*(S_p\wr G)=(CH^*(\Z/p\wr G))^W$ (cf \ref{swan}) if $i$ is a multiple of $(p-1)$, and the image of the map above is indeed contained in the $CH^*BS_p=F_p[v^{p-1}]$-module generated by the $\st x$.
\end{proof}

\end{blank}

\begin{blank}\label{plagiat}{\hl A few properties of $\RE$.} Put $P=CH^*B\Z/p$ and $Q=CH^*BS_p$.

\begin{lem2}[1] If $\{x\}$ is an $\F_p$-basis of $M$, then $\{\st x\}$ is a basis of $\RE M$ as a free $Q$-module.
\end{lem2}

\begin{proof} (Almost word for word from \cite{lz}, proof of 4.2.3, included for the convenience of the reader). Given a relation $\sum \lambda_x \st x=0$ with $\lambda_x\in P$, put $m=\textrm{min}~ |x|$ taken over all $x$ for which $\lambda_x\ne 0$, if there are any. Then project into $P\otimes M^m$, which is a free $P$-module having a basis containing the elements $1\otimes x$ for $x$ of degree $m$. You get, with $m=2k$, the relation  $\sum \lambda_x (-1)^k v^{k(p-1)}\otimes x=0$ which yields $\lambda_x=0$. \end{proof}

\begin{lem2}[2] The functor $\RE$ is exact.
\end{lem2}

\begin{proof} From the previous lemma, we can say that as a vector space, $\RE M$ is $Q\otimes \Phi M$, where here $\Phi M$ means $M$ with all degrees mutiplied by $p$. Result follows.
\end{proof}

\begin{lem2}[3] IF $M'$ is a submodule of $M$, then $(\RE M) \cap (P\otimes M') = \RE M'$.
\end{lem2}

\begin{proof} This follows from the commutative diagram with exact rows:
$$\begin{CD}
0 @>>> \RE M' @>>> \RE M @>>> \RE M/M' @>>> 0 \\
@.       @VVV     @VVV       @VVV           \\
0 @>>> P\otimes M' @>>> P\otimes M @>>> P\otimes M/M' @>>> 0
\end{CD}$$
\end{proof}

Our next lemma will involve the functor $\Phi$ defined in \cite{lionel}, 1.7. When $M\in \UE$, $\Phi M$ is $M$ with all degrees multiplied by $p$, and with an appropriate action of the Steenrod algebra which makes the map $\Phi M \to M, x \mapsto P_0x$ linear over $\A$ (so more formally, $P^i\Phi x = \Phi P^{i/p}x$ if $p|i$ and $0$ otherwise). Whenever $M$ is reduced, we identify $\Phi M$ with the submodule of $M$ comprised of the elements $P_0x$.

\begin{lem2}[4] For any $M\in\UE$, we have
$$(\RE\Phi M) \cap (\Phi P\otimes M) \subset \Phi\RE M$$
\end{lem2}

\begin{proof}(cf \cite{preprint}, 4.7). Let $x\in M$ have degree $2k$; we can write $$\st \Phi x =\sum_{j=0}^{k} (-1)^j v^{jp(p-1)}\otimes P^{(k-j)p} \Phi x$$ because $P^{kp-i}$ can only act non-trivially on $\Phi x$ when $p$ divides $i$ (so we have put $i=pj$). A typical element in $(\RE\Phi M) \cap (\Phi P\otimes M)$ is thus $$y=v^{pm(p-1)}\sum_{j=0}^{k} (-1)^j v^{jp(p-1)}\otimes P^{(k-j)p} \Phi x$$ and this is $\Phi z$ for $$z=v^{m(m-1)}\sum_{j=0}^{k} (-1)^j v^{j(p-1)}\otimes P^{(k-j)} \Phi x$$ from Cartan's formula and 1.7.1* in \cite{lionel}.
\end{proof}

\end{blank}

\begin{propo} If $M\in \UE$ is $\nil$-closed, so is $\RE M$.
\end{propo}

\begin{proof} We prove that if an element $y$ of $\RE M$ is of the form $y=P_ox$ for some $x\in CH^*B\Z/p\otimes M$, then in fact we can choose such an $x$ in $\RE M$. It follows that $P_0$ is injective on the quotient $(CH^*B\Z/p\otimes M)/\RE M$, so that this module is reduced (cf \ref{nildefs}). As $CH^*B\Z/p\otimes M$ is $\nil$-closed from  proposition \ref{nilprop}, 4/, it follows from 2/  of the same propostion that $\RE M$ is $\nil$-closed.

We have
$$(\RE\Phi M) \cap (\Phi P\otimes M) \subset \Phi\RE M$$ and
$$\RE M \cap (P\otimes \Phi M)=\RE\Phi M$$ from (\ref{plagiat}, lemmas 3 \& 4).
Thus
$$\RE M \cap \Phi(P\otimes M) \subset  \Phi\RE M$$ using the fact that on $\UE$, $\Phi$ commutes with tensor products, just like it does on $\U$ at the prime 2, so that $\Phi(P\otimes M)=(\Phi P\otimes M) \cap (P\otimes \Phi M)$.

This was what we wanted.
\end{proof}

\begin{blank} By contrast, proving that $\tau(CH^*BG^{\otimes p})^W$ is $\nil$-closed is straightforward. Put $M=CH^*BG$. The explicit formulae of \ref{bigraded} and \ref{swan} show that there is an exact sequence
$$0\to \tau(M^{\otimes p}) \to (M^{\otimes p})^{\Z/p} \to \Phi M \to 0$$
Since $M$ is reduced, direct computation shows that $\Phi M$ is reduced; the middle term of the sequence is $\nil$-closed by \ref{nilprop},4/ and \ref{invlim}; so $\tau(M^{\otimes p})$ is $\nil$-closed by \ref{nilprop}, 2/; and finally $(\tau(M^{\otimes p}))^W$ is $\nil$-closed by \ref{invlim} again.
\end{blank}

\begin{theom}\label{mainsym} Let $G$ be a group as in the introduction. If $CH^*BG$ is $\nil$-closed, so is $CH^*B(S_p\wr G)$. In particular, $CH^*BS_n$ is $\nil$-closed for any $n$. It follows that $$CH^*BS_n=\invlim CH^*BE$$ where $E$ runs over the elementary abelian subgroups of $S_n$.
\end{theom}

\begin{proof} The only thing to add is that the Sylow subgroup of $S_n$ is contained in a product of iterated wreath products $S_p\wr S_p \wr\cdots\wr S_p$.
\end{proof}

\begin{blank} \label{feedback} It has been known since Quillen's paper \cite{quillen} that the cohomology of $S_n$ is reduced. From \ref{tob}, we deduce $$CH^*BS_n=\tob H^*BS_n$$ It follows that $$\tob H^*BG = \invlim \tob H^*BE$$ which was established in \cite{preprint}, a paper which has been a great source of inspiration.
\end{blank}

\section{Localization of Chevalley groups}

\begin{blank}{\hl Chevalley groups.}\label{mainchev} The point of this section and the next is to study Chevalley groups, ie groups of the form $G(\F_{l^a})$ where $G$ is a connected, reductive, split group scheme over $\Z$, and $l$ is a prime different from $p$. Recall that ``split'' means that $G$ has a maximal torus $T$ which is itself ``split'' in the sense that $T=\mathbb{G}_m^n$ over $\Z$ (or whatever the base is); in this definition $\mathbb{G}_m$ is the multiplicative group scheme $\textrm{Spec}~ \Z[X,X^{-1}]$. The associated Lie group over $\C$ will be denoted $\GC$.

\end{blank}

\begin{blank}{\hl The Weyl group.}\label{weyl} We have to start with some general remarks on group schemes. Recall that normalisers, centralisers, and quotients (at least by diagonalisable subgroups) can be performed over any base scheme, see \cite{sga32}, expos\'e VIII,6/. It will be important to us that these operations commute with base extensions -- for example if $G$ is a group scheme over $S$ and $H$ a subgroup scheme, and if $R\to S$ is any morphism, then $(N_G(H))_R=N_{G_R}(H_R)$ where $X_R=X \times_S R$.

If now $G$ is a reductive group scheme over $S$ with a split maximal torus $T$, recall that $T=C(T)$ (its own centraliser), cf \cite{sga33}, expos\'e XIX, 2.8. We (momentarily) put $W_G(T)=N_G(T)/T$. The following is proved in \cite{sga33}, expos\'e XXII,3/ (in particular proposition 3.4): there exists a finite group $W$ such that $W_G(T)$ is the constant group scheme associated to $W$. This means in particular that for any ring $R$ above $S$ without idempotents other than 0 and 1, we have $W_G(T)(R)=W$. Furthermore, it is also established in {\em loc cit} that $$W\subset \frac{N_G(T)(S)}{T(S)} \subset W_G(T)(S)$$ Using the preceding remark on base extensions, we conclude that the above inclusions are equalities whenever $S$ is replaced by (the spectrum of) a ring without nontrivial idempotents.

Now let $G$ be a Chevalley group as in the previous paragraph. There is thus a finite group $W$ such that 
$$\begin{array}{ccccc}
W & = & \frac{N_G(T)(\Z)}{T(\Z)} & = &  \frac{N_G(T)}{T} (\Z) \\
  & = & \frac{N_G(T)(k)}{T(k)} & = &  \frac{N_G(T)}{T} (k)
\end{array}$$ for any field $k$.

It follows that $W$ acts on $T(k)$. Note that unless $k$ is algebraically closed, there is no reason for $N_G(T)(k)$ to be the normaliser of $T(k)$ in $G(k)$, it may be a strictly smaller subgroup; similarly $T(k)=C(T)(k)$ may be smaller than the centraliser of $T(k)$ in $G(k)$. However we have the following

\begin{lem2} Any automorphism of $T(k)$ induced by conjugation by an element of $G(k)$ may also be realised by an element of $W$. More generally, any isomorphism between two subgroups of $T(k)$ induced by conjugation by an element of $G(k)$ can also be realised by conjugation by an element of $W$.
\end{lem2}

\begin{proof}
This is a well-known argument. So let $K$ denote an algebraic closure of $k$. 
Suppose $x\in G(k)$ and $xAx^{-1}=B$. Let $C$ be the connected component of 1 in the centraliser of $B$ in $G(K)$. We have then
 $T(K)$ and $xT(K)x^{-1}$ as maximal tori in $C$, and therefore they are conjugated by some $c\in C \subset G(K)$.  Consider then $n=c^{-1}x$. It is clear that $n\in N_G(T)(K)$, and can be written $n=wt$ with $w\in N_G(T)(k)$ and $t\in T(K)$, according to the equalities above, valid for any field. Thus the elements $w$ and $n$ induce the same automorphism of $T(K)$, and induce the same isomorphism between $A$ and $B$ as $x$ does. \end{proof}

\end{blank}

\begin{blank}{\hl Notations.} Throughout the rest of the paper, $G$ will be a Chevalley group, $T$ a split maximal torus, $N_T$ will be short for $N_G(T)$, $W$ will be the finite group just introduced (to be refered to as the Weyl group), and finally $k$ will be a finite field of characteristic different from $p$ but containing the $p$-th roots of unity. It is assumed that $H^*(B\GC,\Z)$ has no $p$-torsion, and moreover that $p$ is odd. Some of this will be repeated for emphasis.
\end{blank}

\begin{blank}{\hl Localization of Chevalley groups.}\label{pedroborel} The elementary abelian subgroups of Lie groups have been studied extensively by Borel \cite{borel}, among others. He proves the following amazing theorem:

\begin{thm2} Let $G$ be a compact connected Lie group. Then the following three conditions are equivalent:
\begin{enumerate}
\item $H^*(G,\Z)$ has no $p$-torsion,
\item $H^*(BG,\Z)$ has no $p$-torsion,
\item any elementary abelian subgroup of $G$ is contained in a maximal torus.
\end{enumerate}
\end{thm2}

The reader might find it amusing to see how the implication 2/ $\Rightarrow$ 3/ can be proved very quickly using some recent results in homotopy theory: namely, if $E$ is elementary abelian, then conjugacy classes of maps $E\to G$ are precisely the same as maps $H^*BG \to H^*BE$ of unstable algebras (\cite{lannes}, corollaire 3.1.4); the assumption on $p$ implies, classically, that $H^*BG\to H^*BT$ is injective, where $T$ is any maximal torus, and thus what we need to prove is that we can extend the map $H^*BG\to H^*BE$ induced by the inclusion of $E$ into $G$ to a map $H^*BT\to H^*BE$ of unstable algebras. In turn, this follows from the ``non-linear injectivity of $H^*BE$'', \cite{lionel}, 3.8.7.

We prove now a variation on Borel's theorem. 

\begin{prop2} Let $G$ be a Chevalley group with split maximal torus $T$, let $p$ be a prime number such that $H^*(B\GC,\Z)$ has no $p$-torsion, and let $k$ be a finite field of characteristic $\ne p$ which contains the $p$-th roots of unity. Then any elementary abelian $p$-subgroup of $G(k)$ is conjugated to a subgroup of $T(k)$.
\end{prop2}

\begin{proof} It is proved in \cite{pedro} that for such $k$, we have 
$$\MUP{BG(k)}=H^*B\GC$$ for any $G$ (hence for $T$ as well). By choice of $p$, the map $H^*B\GC \to H^*BT_\C$ is injective, so that
$$\MUP{BG(k)}\to\MUP{BT(k)}$$ is injective. Combining this with (\ref{cobordism}, proposition) tells us that the natural map
$$\underset{\CG}{\invlim} \MUP{BE} \rightarrow \underset{\CdG}{\invlim} \MUP{BE}$$ is an F-monomorphism, where $\CdG$ is the subcategory of $\CG$ (definition in section \ref{section:nil}) consisting of those $E$ which are conjugated to a subgroup of $T(k)$. (In fact this is simply a monomorphism, as both its source and target are $\nil$-closed modules in $\UE$.)

Suppose that there is a $V\in\CG$ which is not in $\CdG$. Choose a $V$ maximal with respect to this property, and note that $V$ is then maximal in $\CG$. We construct an element $x$ in $\MUP{BV}$ which restricts to $0$ in any proper subgroup of $V$: for this, choose for each such subgroup $E$ a non-zero element $x_E$ in $\MUP{BV}$ which restricts to $0$ in $\MUP{BE}$, for example the first Chern class of a non-trivial $1$-dimensional representation of $V$ which is trivial on $E$; and then take $x$ to be the product of all the different $x_E$'s. By symmetry, $x$ is invariant under the action of the normaliser of $V$. (In fact, in this way we end up taking $x$ to be the product of {\em all} non-zero elements of degree 2, but I prefer to phrase it as I have.)

Given this, define $\alpha(E)=0$ if $E\in\CG$ is not conjugated to $V$ and $\alpha(E)=x$ otherwise, with an obvious abuse of notation. By maximality of $V$ and choice of $x$, this defines an element $\alpha$ in the inverse limit on the left hand side above (we need such an $x$ because there could be a group $W$ which is not conjugated to $V$ but having a subgroup $E$ conjugated to a subgroup of $V$).

This is a contradiction, as $\alpha$ is non-zero (ie non-nilpotent) but it lies in the kernel of the (F-) monomorphism above. Hence $\CG=\CdG$.\end{proof}

\begin{coro2} The localization of $CH^*BG(k)$ away from $\nil$ is $$L(CH^*BG(k))=(CH^*BT(k))^W$$ 
\end{coro2}

\end{blank}

\begin{rmk2}[1] Here's a variant of the proposition. Keep the same hypotheses but do not assume that $k$ contains the $p$-th roots of unity; then each elementary abelian $p$-subgroup of $G(k)$ is contained in a maximal torus, possibly not split. To see this, let $E$ be such a subgroup, and let $K$ be an algebraic closure of $k$. Let $Z(E,G)$ be the centraliser of $E$ in $G(K)$, and let $Z^0$ be the connected component of 1 -- a reductive group. From the result obtained when $k$ is big enough, we deduce that $E$ is contained in $Z^0$. Being central, it is contained in any maximal torus of $Z^0$, and there is one which is defined over $k$. \end{rmk2}

\begin{rmk2}[2] Our result can be recovered from a theorem of Steinberg, as follows. Let $E$ be an elementary $p$-subgroup of $G(k)$ and let $K$ be an algebraic closure of $k$. Then theorem 2.28 in \cite{steinberg} says that $E$ is toral over $K$, ie there is a $g\in G(K)$ such that $gEg^{-1}\subset T(K)$ where $T$ is the fixed, split maximal torus. Considering the assumption on $k$, $gEg^{-1}$ is in fact contained in $T(k)$. Now consider the set $X$ of $x\in G(K)$ such that $xex^{-1}=geg^{-1}$ for all $e$ in $E$. It is a principal homogeneous space under $Z$, the centraliser of $E$ in $G(K)$. This $Z$ is defined over $k$ and connected (theorem 2.28 in {\em loc cit} again), and by Lang, it has a $k$-rational point. In other words there is an $x\in G(k)$ such that $xEx^{-1}$ is contained in $T(k)$, which was what we wanted.
\end{rmk2}

\begin{note2} The two remarks above are due to J.P. Serre.
\end{note2}

\begin{blank}  We will call condition (1), resp.  (2), the injectivity, resp. surjectivity, of $$CH^*BG(k)\to (CH^*BT(k))^W$$ 
We will see shortly that (1) or (2) holds for $G$ if it holds for its universal cover $\tG$, with a partial converse. 

It might be enlightening, to start with, to recall what is known about the cobordism of Chevalley groups.
\end{blank}

\begin{blank}{\hl Cobordism.}\label{munilclosed} The cobordism of Chevalley groups is completely described in \cite{pedro}; in particular, the result already quoted says
$$\MUP{BG(k)}=H^*B\GC$$ (a natural isomorphism with respect to maps of group schemes over $\Z$.) It is also established that the latter rings inject into $H^*BG(k)$, and are polynomial.

We deduce first of all that $\MUP{BG(k)}$ can be viewed as an element in $\UE$. If $p$ is odd, we can say a bit more: a theorem of Dwyer-Miller-Wilkerson (\cite{dwyer}, 2.11) implies then that $$\MUP{BG(k)}=(\MUP{BT(k)})^W$$ (in fact equality holds if and only if the left hand side is polynomial, cf \cite{moller}, lemma 7.1.) In other words, $\MUP{BG(k)}$ is $\nil$-closed in $\UE$.
\end{blank}

\begin{blank}\label{covers}{\hl Covers.} Assume that $G$ is semi-simple and consider now the universal cover $\tG$ of $G$, ie the unique Chevalley group such that $\tGC$ is the universal cover of $\GC$. The centre of a semi-simple Chevalley group is finite and contained in any maximal torus, hence is a product of groups $\mu_a$ ($a$-th roots of unity). These groups have the peculiarity that all subgroups of $\mu_a(\C)$ are defined over $\Z$, and the same can be said thus of the subgroups of the centre $\mathcal Z(\C)$ of $\tGC$. In particular $\pi_1(\GC)$ can be seen as a finite and central subgroup $\pi$ of $\tG$, defined over $\Z$. Now, the groups $\tG/\pi$ and $G$ agree over $\C$, and hence agree over $\Z$, by the structure theorem for Chevalley groups, cf \cite{sga33}, expos\'e XXV, th\'eor\`eme 1.1.

The map $\tG (k) \to G(k)$, where as always $k$ is finite, may not be as nice as one would expect it to be (eg it is almost never surjective). Its kernel, at least, is certainly a subgroup of $\pi(k)$, so it is central. The order of $\mu_a(k)$ divides that of $\mu_a(\C)$, and it follows that $|\pi(k)|$ divides $|\pi(\C)|=|\pi_1(\GC)|$ (this is a result that holds for more general groups, but it is trivial to check it directly here). Now let $Tors(G)$ denote the set of prime numbers $l$ such that $H^*(B\GC,\Z)$ has $l$-torsion. All prime numbers dividing the order of $\pi_1(\GC)$ are in $Tors(G)$, in fact $Tors(G)=Tors(\tG)\cup \{l: l~ \textrm{divides}~ |\pi_1(\GC)|\}$: see \cite{serre}, 1.3.1 to 1.3.4, and the references there, in particular \cite{borel}. So our standing assumption that $H^*(B\GC,\Z)$ have no $p$-torsion implies that the order of $\pi(k)$, and thus that of the kernel of $\tG (k) \to G(k)$, is prime to $p$.

We also have to note that $|\tG (k)|=|G(k)|$: this is very specific to finite fields, and is proved in \cite{boreltxt}, proposition 16.8. It follows that the image of $\tG(k)$ in $G(k)$ has index prime to $p$. 

We note:

\begin{lem2} Let $\Gamma$ be given as a quotient $\tilde\Gamma/C$ where $C$ is a central subgroup of the finite group $\tilde\Gamma$. Assume that the order of $C$ is prime to $p$. Then $H^*\Gamma=H^*\tilde\Gamma$ and $CH^*B\Gamma=CH^*B\tilde\Gamma$.
\end{lem2}

\begin{proof} A spectral sequence argument gives the result immediately for cohomology. The result for Chow rings follows from remark \ref{martino}, but of course it is preferable to find an elementary argument (even if we only use the easier part of \ref{martino}). We let $\tilde S$ and $S$ denote Sylow subgroups of $\tilde\Gamma$ and $\Gamma$ respectively. Clearly $\tilde S$ maps isomorphically to (a conjugate of) $S$ under the quotient map; we now use proposition \ref{stable}. To prove that the stable subgring for $S$ is sent isomorphically onto the stable subring for $\tilde S$ under the isomorphism $CH^*BS \to CH^*B\tilde S$, things come down to checking the following: if $\tilde x\in\tilde S$ and if $\tilde g\in\tG(k)$, then $\tilde g\tilde x\tilde g^{-1}$ is in $\tilde S$ if and only if $gxg^{-1}$ is in $S$, where $x$ and $g$ correspond to $\tilde x$ and $\tilde g$ under the projection $\tilde\Gamma\to \Gamma$. To see the non-trivial part of this, assume that $gxg^{-1}\in S$ and write $\tilde g\tilde x\tilde g^{-1}=c\cdot s$ with $c\in\pi(k)$ and $s\in\tilde S$. Writing $o(\sigma)$ for the order of an element $\sigma$, we have $o(c\cdot s)=o(c)o(s)$ because $c$ and $s$ commute and have coprime orders. But $o(c\cdot s)$ is a power of $p$, as this element is conjugate to $\tilde x$. Therefore $o(c)=1$ and $c=1$.\end{proof}

Applying this with $\Gamma$ taken to be the image of $\tG(k)$ in $G(k)$ shows that $CH^*BG(k)$ is a direct summand in $CH^*B\tG(k)$ (in $\UE$). So if $CH^*B\tG(k)$ is reduced or $\nil$-closed, so is $CH^*BG(k)$. In fact, we can say something more precise. 

Letting $\Gamma$ be as above, we claim that $G(k)$ is the product of $T(k)$ and $\Gamma$ (although these subgroups have a non-trivial intersection, of course). To see this, let $K$ be an algebraic closure of $k$; we may assume to have chosen a maximal torus $\tilde T$ in $\tG$ such that $p(\tilde T(K))=T(K)$, where $p: \tG\to G$ is the natural map. Since $\tilde T(k)$ contains the kernel $C$ of the map $\tG(k)\to G(k)$, it follows easily that the image $E$ of $\tilde T(k)$ in $T(k)$ is all of $\Gamma\cap T(k)$. Moreover, since $\tilde T(k)$ and $T(k)$ have the same order (by the result already quoted, or directly), it follows that $[T(k):E]=[G(k):\Gamma]$ (in turn these indices equal the order of $C$). Consequently, the number of products $\gamma\cdot t$ with $\gamma\in\Gamma$ and $t\in T(k)$ is 
$$\frac{|\Gamma|\cdot |T(k)|}{|\Gamma\cap T(k)|}=|\Gamma|\cdot [T(k):E]=|\Gamma|\cdot [G(k):\Gamma]=|G(k)|$$

Thus $G(k)=\Gamma\cdot T(k)$, as claimed. We apply now the double coset formula (see the appendix), which gives here:
$$i^*_{T(k)\to G(k)}\circ i_*^{\Gamma \to G(k)}=i_*^{E\to T(k)}\circ i^*_{E\to\Gamma}$$
Since $E$ has index prime to $p$ in the abelian group $T(k)$, it is clear that $i_*^{E\to T(k)}$ is an isomorphism. The transfer $i_*^{\Gamma \to G(k)}$ is surjective for index reasons. Thus we see from the formula above that the surjectivity of $CH^*BG(k)\to (CH^*BT(k))^W$ is equivalent to that of $CH^*B\Gamma\to (CH^*BE)^W$. From the last lemma applied to $\Gamma$ and $E$, this is also equivalent to $\tG$ satisfying (2). Summarizing:

\begin{prop2} Let $\tG$ be the universal cover of the semi-simple Chevalley group $G$. If  $\tG$ satisfies (1), resp. (2), then $G$ satisfies (1), resp. (2). Moreover if $G$ satisfies (2), then so does $\tG$.
\end{prop2}

\begin{rmk2} It is clear from the discussion above that the proposition holds a bit more generally: if $\tG\to G$ is an isogeny (over $\Z$) whose kernel $C$ is central and has $|C(k)|$ prime to $p$, then the conclusion of the proposition holds.
\end{rmk2}

\begin{rmk2}[on notation] The simply-connected group $Spin_n$ has a central subgroup $\mu_2$, and the quotient $Spin_n/\mu_2$ is the orthogonal group $SO_n$. Note that the group we call $SO_n(k)$ is usually denoted $SO^+_n(k)$, while $SO_{2n}(\R)$ sometimes goes under the name $SO(n,n)$ (and $SO_{2n+1}(\R)$ can be $SO(n+1,n)$). The familiar group of orthogonal motions of $\R^n$ with determinant 1, normally denoted $SO(n)$, is the compact form of the group we consider.
\end{rmk2}

\section{Chevalley groups: examples}

\begin{blank}\label{pnotdivw}{\hl When $p$ does not divide $|W|$.} We assume that $G$ is semi-simple and as usual that the characteristic of $k$ is not $p$. Under these hypotheses, a theorem of Springer-Steinberg (\cite{springer} 5.19) asserts that a $p$-Sylow $S$ of $G(k)$ normalises a maximal torus $T'$ defined over $k$. This $T'$ does not have to be split; however if we suppose further that $p$ does not divide the order of the Weyl group of $G$, it follows that $S$ is abelian. If we denote its normaliser in $G(k)$ by $N_S$, we have by Swan's lemma (\ref{swanchow}) $$CH^*BG(k)=(CH^*BS)^{N_S}$$ 
Now, let $S'$ be the maximal $p$-elementary subgroup of $S$. We have (cf \ref{intronotations}) $CH^*BS=CH^*BS'$. Let $N_{S'}$ denote the normaliser of $S'$ in $G(k)$. Clearly $N_S\subset N_{S'}$ (since $S'$ is precisely the subgroup of $S$ of elements of order $p$, it is preserved by any automorphism of $S$) and thus $(CH^*BS')^{N_{S'}}\subset (CH^*BS')^{N_S}$. However since there are isomorphisms
$$CH^*BG(k)\to (CH^*BS)^{N_S} \to (CH^*BS')^{N_S}$$ we obtain the reverse inclusion, and $CH^*BG(k)=(CH^*BS')^{N_{S'}}$.

Now, $S'$ being a $p$-elementary subgroup of $S$ of maximal rank, it is also of maximal rank in $G(k)$; since (by \ref{pedroborel}) we know that $S'$ is (conjugated to) a subgroup of $T(k)$ where $T$ is our fixed, split maximal torus, we deduce that $S'$ is precisely the maximal elementary abelian subgroup of $T(k)$ and that $CH^*BT(k)=CH^*BS'$ (\ref{intronotations} again).

Finally, the group of automorphisms of $CH^*BS'$ induced by $N_{S'}$ is the same as that induced by $N_T(k)$ (or $W$), by (\ref{weyl}, lemma). Thus:

\begin{prop2} Suppose $G$ is semi-simple. If $p$ does not divide $|W|$, then $CH^*BG(k)$ is $\nil$-closed, ie $$CH^*BG(k)=(CH^*BT(k))^W$$
\end{prop2}

\begin{ex2} Consider the five exceptional Lie groups $G_2$, $F_4$, $E_6$, $E_7$ and $E_8$. Their Weyl groups have orders  $12$, $2^7.3^2$, $2^7.3^4.5$, $2^{10}.3^4.5.7$ and $2^{14}.3^5.5^2.7$ respectively. The primes which give torsion in their integral cohomology are $\{2\}$ for $G_2$, $\{2, 3\}$ for $F_4$, $E_6$ and $E_7$, and $\{2, 3, 5\}$ for $E_8$ (see \cite{serre}). In other words, we can say that $CH^*BG(k)$ is always $\nil$-closed when $G$ is exceptional and does not have $p$-torsion, except possibly in finitely many cases, namely when $(G,p)$ is either $(G_2, 3)$, $(E_6, 5)$, $(E_7, 5)$, $(E_7, 7)$, or $(E_8,7)$. Restricting attention to $p>7$ rules out these cases.
\end{ex2}

\end{blank}

\begin{blank}\label{chern}{\hl Chern classes.} Another case of interest is that of a group $G$ such that $H^*B\GC$ is generated by Chern classes of (finitely many) representations of $\GC$. In this case we have:

\begin{prop2} Suppose that $H^*B\GC$ is generated by Chern classes, and let $p$ be odd. Then the map $$CH^*BG(k)\to \MUP{BG(k)}$$ is surjective, as is  $$CH^*BG(k)\to (CH^*BT(k))^W$$
\end{prop2}

\begin{proof} By assumption there is a surjective map
$$\bigotimes_i H^*BGL(n_i,\C) \to H^*B\GC$$  Such representations of $\GC$ must be defined over $\Z$ (essentially because irreducible representations of $G$ over $\Z$ are classified by their highest weights, just like those of $\GC$ over $\C$; and as $G$ is split, the weights are the same over $\Z$ or $\C$. Alternatively, appeal to \cite{sga33}, XXV, 1.1.).  Consider the following commutative diagram:
{\tiny
$$\begin{CD}
\bigotimes CH^*BGL(n_i,k) @>>> CH^*BG(k) @>c>> (CH^*BT(k))^W \\
    @VaV{\approx}V               @VbVV           @VV{\approx}V \\
\bigotimes \MUP{BGL(n_i,k)} @>d>> \MUP{BG(k)} @>e>> (\MUP{BT(k)})^W 
\end{CD}$$}
The map $a$ is an isomorphism by the computation of $CH^*BGL(n,k)$ that was made in \cite{pedro}. Moreover $d$ is surjective by the remark just made and the results of {\em loc cit}. It follows that $b$ is surjective. Finally,  $e$ is an isomorphism by (\ref{munilclosed}), so $c$ is surjective. \end{proof}

\begin{ex2} Groups such that $H^*B\GC$ is generated by Chern classes include $SL_n$, $Sp_n$, and $SO_{2n+1}$ for odd $p$'s. Here again we point out that $SO_{2n+1}$ has 2-torsion in its integral cohomology, so the restriction that $p$ be odd in the above proposition is not important anyway; however for $SL_n$ and $Sp_n$, we are leaving out cases where the conclusion of the proposition might still be true. 
\end{ex2}

\end{blank}

\begin{blank}{\hl The orthogonal groups.} We have already established that, with some restrictions on $p$, all simple, simply connected groups satisfy condition (2), with the notable exception of $Spin_{2n}$. Replacing this group by $SO_{2n}$ as we may, we can use an {\em ad hoc} argument.

\newcommand{\sok}{SO_{2n}(k)}
\newcommand{\soc}{SO_{2n}{\C}}

 We refer back to \cite{pedro} where it is proved that a certain continuous map $B\sok\to B\soc$ induces isomorphisms
$$\begin{CD}
\MUP{B\soc} @>{\approx}>> \MUP{B\sok} @>{\approx}>> (\MUP{BT(k)})^W
\end{CD}$$

Each of these rings is a polynomial ring over $\F_p$ on variables $c_{2i}$ ($1\le i \le n-1$) and $e$. By definition then, these classes are the Pontryagin classes and the Euler class of the oriented bundle $E\to B\sok$ induced by the classifying map $B\sok\to B\soc$. Now if this bundle were algebraic, it would have Pontryagin and Euler classes in $CH^*B\sok$ which would restrict to the elements with the same name in $CH^*BT(k)$, proving that (2) holds (recall that $CH^*BT(k)=\MUP{BT(k)}$, of course, and that the cobordim modules are $\nil$-closed, cf \ref{munilclosed}). However $E$ need not be algebraic.

The Pontryagin classes are Chern classes and one can argue as in \ref{chern} to prove that they exist in $CH^*B\sok$. We need more work for the Euler class.

So let $S$ be a $p$-Sylow of $\sok$. It is proved in \cite{DZ} that any map $BP\to BG$, where $p$ is a finite $p$-group and $G$ a compact Lie group, comes from a homomorphism $P\to G$, up to homotopy. Therefore the composition $BS\to B\sok\to B\soc$ is homotopic to a map coming from a homomorphism $S\to \soc$. As $S$ is finite, this is automatically a map of algebraic groups, and the map $BS\to B\soc$ is algebraic. Using this map to pull back the universal vector bundle over $B\soc$, we get an algebraic vector bundle over $BS$ carrying a quadratic form and which is equivalent to $E$ (or rather, to the pull back of $E$ over $BS$). By \cite{edidin}, this bundle has an Euler class in the Chow ring, ie there is a class in $CH^*BS$ mapping under the cycle map to the Euler class in $\MUP{BS}$ (up to scalar multiplication by a power of 2, but we are still assuming that $p$ is odd).  Call this $x\in CH^*BS$, put $e'=i_*^{S\to\sok}(x)$, and finally let $e=i^*_{T(k)\to\sok}(e')$. To prove that this is the Euler class defined above (and justify the notation), it is enough to check this in $\MUP{BT(k)}$, that is, after applying the cycle map. But in cobordism we have Euler classes for topological bundles, and $x=i^*_{S\to\sok}(y)$ for some $y\in \MUP{B\sok}$, so that $e'=y$ up to a nonzero scalar (we abuse the notations $x$ and $e'$ here). It follows that $e$ is indeed the Euler class that we want.

\end{blank}

\begin{blank}{\hl Condition (1) and wreath products.} \label{wreath} The Weyl groups of the simple, simply-connected Lie groups have a lot in common. Namely, for the classical groups at least, $W$ always contains a copy of the symmetric group acting on the torus by permuting the eigenvalues. This often gives $N_T(k)$ a tractable group structure. The following lemma is trivial, but it is remarkable to note how often it applies.

\begin{lem2}
 Suppose that $G(k)$ contains a subgroup of index prime to $p$ which is of the form $S_r\ltimes T(k)=S_r \wr k^*$ (here $r$ is the rank of $G$). Then the restriction map:
$$CH^*BG(k) \to CH^*BT(k)$$ is injective. 
\end{lem2}

\begin{proof}(Compare with \ref{mainsym}.) The restriction to $S_r\wr k^*$ is injective for index reasons. If follows from \cite{pedro}, lemma 4.4, that $S_r\wr k^*$ possesses a collection of abelian $p$-groups which detect the Chow ring -- in the sense that if $x\in CH^*B(S_r\wr k^*)$ restricts to $0$ in each of these subgroups then $x=0$. Thus the Chow ring is also detected, in the same sense, by the family of elementary abelian p-subgroups. As observed above, these are conjugated to subgroups of $T(k)$.
\end{proof}

\begin{rmk2} The reader will have noticed that an easy induction based on \ref{mainsym} will suffice to prove that $S_r \wr k^*$ has a reduced Chow ring. However, the proof in \cite{pedro} and a number of results in the same paper yield the slightly stronger statement that this group possesses a subgroup $H$ of index prime to $p$ whose Chow ring is reduced and which satisfies a certain condition (called (*) there) meaning in fact the hypotheses of \ref{intro} together with the property that the cycle map to cohomology is injective. In turn, two such groups $H_1$ and $H_2$ have a K\"unneth formula, in the sense that $$CH^*BH_1\times BH_2=CH^*BH_1\otimes CH^*BH_2$$ which is not true for the Chow ring of a product of varieties in general, but see {\em loc cit} and \cite{totaro} in this case. 

It follows that if the lemma applies for two Chevalley groups $G_1$ and $G_2$, then their product also has a reduced Chow ring.
\end{rmk2}

\begin{ex2}[1] It is easily seen that the lemma applies, for odd $p$, for the groups $Sp_n$ and $Spin_n$.
\end{ex2}

\begin{ex2}[2] In some cases we can also deal with $SL_n$. Here $N_T(k)$ is of the form $T(k)\rtimes S_n$, but the action of $S_n$ is trickier: if we see $T(k)$ as $\{(x_1,\cdots,x_n)\in (k^*)^n : x_1\cdots x_n=1\}$ (that is, as a subgroup of the maximal torus for $GL_n$), then $S_n$ acts by permutation of the $x_i$'s. However, when $n$ is prime to $p$, then we can have a look at the subgroup $T\rtimes S_{n-1}$ where $S_{n-1}$ sits in $S_n$ as the subgroup which fixes, say, the $n$-th coordinate. This subgroup is indeed a wreath product $S_{n-1}\wr k^*$ of index prime to $p$ in $N_T(k)$, and when $p$ is also assumed to be odd then $N_T(k)$ has order prime to $p$ in $SL_n(k)$.
\end{ex2}

\end{blank}

\end{blank}

\begin{blank}{\hl Summary.} Let us put together the information gathered so far. First of all, there is the

\begin{thm2}\label{main} Let $G$ be a Chevalley group, let $p$ be a prime number, and let $k$ be a finite field of characteristic $\ne p$ containing the $p$-th roots of unity. Assume that $G$ has no $p$-torsion. Then $CH^*BG(k)$ is $\nil$-closed in $\UE$ if $G$ is locally isomorphic to a product of the following groups:

\begin{itemize}
\item{Type $A_{n-1}$:} $GL_n$ for $p>2$;  $SL_n$ if $n$ is prime to $p$ and $p>2$;
\item{Type $B_n$ or $D_n$:} $Spin_n$ for $p>2$;
\item{Type $C_n$:} $Sp_n$ for $p>2$;
\item{Exceptional type:} $G_2$ for $p>3$; $F_4$ for $p>3$; $E_6$ for $p>5$; $E_7$ for $p>7$; $E_8$ for $p>7$.
\end{itemize}
\end{thm2}

Here ``locally isomorphic'' refers to the situation in \ref{covers} (see the remark following the proposition there: this is why it makes sense to include $GL_n$, which is not simply-connected, in the list).

Our main omission with simply-connected groups is $SL_n$ when $p$ divides $n$. If one could prove the result in this case as well, we would be able to say that $CH^*BG(k)$ is $\nil$-closed for {\em any} semi-simple group (and say, for $p>7$).

Nevertheless, we have 

\begin{prop2}\label{image} Let $G$ be any semi-simple Chevalley group having no $p$-torsion, and let $k$ be any finite field of characteristic $\ne p$. Let $K$ be the finite field obtained from $k$ by adjoining the $p$-roots of unity and let $r=[K:k]$.

Then the map $$CH^*BG(k)\to \MUP{BG(k)}$$ is surjective. Moreover for a certain Brauer lift $BG(k)\to B\GC$, the image of $H^*B\GC \to H^*BG(k)$ coincides with the image of the Chow ring under the cycle map. Writing $H^*B\GC=\F_p[s_1,\cdots,s_n]$, this image is $$\F_p[s_i: 2r ~divides~ |s_i|].$$
\end{prop2}

\begin{proof} The point is that we have established (2) for any simply-connected group, hence for any semi-simple group by \ref{covers}. The first statement follows immediately from this and \ref{munilclosed}, at least for $K$; but the map $\MUP{BG(K)}\to \MUP{BG(k)}$ is always surjective as we showed in \cite{pedro}. 

Everything else in the proposition was proved in {\em loc cit} for cobordism, with the exception of the injectivity of $$\MUP{BG(k)}\to H^*BG(k)$$ which we need to describe the image and which is only stated there for $K$. The more general statement is obvious though, since $\MUP{BG(k)}$ is a polynomial ring for any $k$ and since the cycle map to cohomology in an $F$-monomorphism. 
\end{proof}

\end{blank}

{\bf Acknowledgements.} This work will be part of my PhD thesis and I would like to thank Burt Totaro for his help as my supervisor. I am also indebted to J-P. Serre for pointing out alternative approaches to the study of elementary subgroups of Chevalley groups (see the two remarks in \ref{pedroborel}).

\appendix

\section{Appendix: the double coset formula}

\begin{blank} The double coset formula is a useful equality that shows up whenever a functor $h^*(-)$ from groups to rings presents itself equipped with transfers. Of course, examples for $h^*(-)$ include $CH^*B-$, $\MUP{B-}$ and $H^*(-)$; example for transfers include not only the standard one, but also Even's norm used by Quillen in \cite{quilewen} and by Yagita in \cite{yagita}.

The lemma that follows implies the double coset formula (to be stated shortly) in all cases of interest. We submit a proof because it is perhaps not well-known that it can be stated in the category of algebraic varieties (as opposed to topological spaces), which is what we need when working with Totaro's definition of a classifying space (\cite{totaro}). In turn, this will fill a very minor gap in Yagita's proof alluded to in section \ref{section:nil}, which implicitly uses the double coset formula for Chow rings and for the Even's multiplicative norm.
\end{blank}

\begin{blank} Let $G$ be a finite group, and let $K$, $H$ be subgroups. Choose representative elements $\sigma_i$ of the $K$-$H$-double cosets, that is, such that G is the disjoint union of the sets $K\sigma_i H$. The intersection $K\cap \sigma_i H \sigma_i^{-1}$ will be written $L_i$.

We let $U$ be a (Zariski) open set in a vector space, on which $G$ acts freely. For each subgroup $A$ of $G$, the variety $U/A$ is Totaro's approximation to $BA$, and this is what the notation will mean in what follows.  

\begin{lem2}
The following square is a pull-back:
$$\begin{CD}
\coprod_i BL_i @>>>  BH \\
@VVV              @VVV \\
BK @>>>  BG
\end{CD}$$

\end{lem2}

\begin{proof} All maps should have an obvious definition. Let $X=\{(x,y)\in BK\times BH : p_1(x)=p_2(y)\}$ where $p_1$ and $p_2$ are the obvious maps to $BK$ and $BH$ respectively. Of course $X=BK\times_{BG} BH$. We first define a map
$BL_i=U/L_i \longrightarrow X\subset U/K \times U/H$
by $\bar u\mapsto (\bar u, \overline{\sigma_i^{-1}\cdot u})$. This is well-defined and takes values in $X$; moreover one checks easily that it gives a bijection $\coprod_i BL_i \to X$. As we work over a field of characteristic $0$, and as all varieties encountered are normal (in fact, smooth -- for example $X$ is \'etale over $BK$ or $BH$), this is an isomorphism of algebraic varieties.
\end{proof}

\end{blank}

\begin{blank}{\hl The formula} As we have just pointed out, there are several versions of the double coset formula, and generalisations to any algebraic group over any field, etc. In all cases it reads:
$$i_{K\to G}^*\circ i_*^{H\to G} = \sum_i i^{L_i\to K}_*\circ i_{L_i \to \sigma_i H \sigma_i^{-1}}^* \circ c_{\sigma_i^{-1}}$$ 
One deduces it easily from the last lemma (cf prop 1.7 in \cite{fulton} for Chow rings).

We now list a few consequences. Proofs are standard and will be omitted.
\end{blank}

\begin{propo}\label{stable}
Let $G$ be a finite group, let $G_p$ be a subgroup of index prime to $p$. Then $CH^k(G)$ restricts isomorphically onto the ``stable'' subgroup of $CH^k(G_p)$. In other words, defining a map $CH^k(G_p)\to CH^k(G_p\cap G_p^g)$ by
$$\Psi_g(x)=i^*_{G_p\cap G_p^g \to G_p}(x) - i^*_{G_p\cap G_p^g \to G_p^g}\circ c_{g^{-1}}^*(x)$$ then the following sequence is exact:

$$0\to CH^k(G)\to CH^k(G_p)\to \prod_{g\in G} CH^k(G_p\cap G_p^g)$$
\end{propo}

\begin{rmk}\label{martino} Let $G$ and $G'$ be finite groups, and let $S$ and $S'$ be Sylows. The Martino-Priddy conjecture, recently proved by Bob Oliver (\cite{bob}), asserts that whenever there is a homotopy equivalence $\p{BG}\approx \p{BG'}$, then there is an isomorphism of groups $S\approx S'$ which is ``fusion preserving'', ie it respects the conjugacies in the strongest possible sense. It follows then from the last proposition that there is an isomorphism $CH^*BG\approx CH^*BG'$. Note that any map $G\to G'$ inducing an isomophism on cohomology yields a homotopy equivalence between the $p$-completed classifying spaces, and thus {\em any homomorphism of groups which induces an isomorphism on cohomology also implies the existence of an isomorphism between the Chow rings} (or cobordism rings, or any ring coming from a theory with transfers for that matter).

However we shall refrain from using this result, as the proof of the Martino-Priddy ``conjecture'' is very long and complicated and proceeds by a case-by-case check through the classification of finite simple groups.
\end{rmk}

\begin{propo}[Swan's lemma] \label{swanchow}Suppose that the $p$-Sylow subgroup $G_p$ is abelian. Then if $N_p$ is its normalizer, we have $$CH^*(G)=CH^*(G_p)^{N_p}$$
\end{propo}

\bibliography{myrefs}
\bibliographystyle{siam}

\end{document}